\documentclass[oneside]{amsart}
\usepackage{amsfonts}
\usepackage{amssymb}
\usepackage{amsxtra}
\usepackage{amstext}
\usepackage[english]{babel}

\newcommand\bra{\left<}
\newcommand\ket{\right>}
\newcommand{\R}{\mathbb{R}}
\newcommand{\C}{\mathbb{C}}

\newcommand{\N}{\mathbb{N}}

\newcommand{\Z}{\mathbb{Z}}

\newcommand{\kO}{\mathcal{O}}

\newcommand{\kX}{\mathcal{X}}

\newcommand{\ka}{\mathfrak{a}}

\newtheorem {lem} {Lemma} [section]
\newtheorem {prop} {Proposition} [section]
\newtheorem {theo} {Theorem} [section]

\newtheorem {rem} {Remark} [section]
\newtheorem {note} {Note} [section]

\newcommand\la{\lambda}
\newcommand\teta{\theta}

\title[Heat kernel and Green function estimates]
      {Heat kernel and Green function\\ estimates on affine buildings of type $\tilde{A}_r$}
\author{Jean-Philippe Anker}
\author{Bruno Schapira}
\author{Bartosz Trojan}

\begin{document}

\begin{abstract} We obtain a global estimate of the transition density $p^n(0,x)$ associated to a nearest neighbor random
walk, called here "simple", on affine buildings of type
$\widetilde{A}_r$. Then we deduce a global estimate of the Green
function. This is the analogue of a result on Riemannian symmetric
spaces of the noncompact type.
\end{abstract}

\maketitle
\let\languagename\relax


\bigskip

\noindent \textbf{A.M.S. classification:} Primary: 42C05; 51E24;
60B15; 60J10.\\ Secondary: 31C12; 31C35; 33C52; 33E50; 60J45;
60J50.

\section{Introduction}
This work is meant as a first attempt to understand the full
behavior of random walks on affine buildings of higher rank. Such
a study was carried out by Lalley (see \cite{L1}, \cite{L2} and
the report in \cite{W}) for rather general random walks on rank
one buildings i.e.~homogeneous trees, and by the first author (in
the joint works \cite{AJ}, \cite{AO1}, \cite{AO2}) for the heat
diffusion on a general Riemannian symmetric space of the
noncompact type, which is a continuous counterpart of the present
discrete setting. Apart from its own interest, this information is
pivotal for further study. For instance, in potential theory, it
is used to estimate the Green's function and to describe the
Martin boundary. It will also be used in \cite{S2} to study the
asymptotic behavior of normalized bridges.

Our paper deals with a simple higher rank case. We consider
buildings of type $\widetilde{A}_r$, which are known to be most
simple among affine buildings, and a particular random walk to the
nearest neighbors, that we shall call ``simple''. This random
walk, actually its Fourier transform, satisfies a ``magic''
combinatorial formula, which is technically very helpful. Our main
result is a global upper bound for the transition density
$p^n(x,y)$, which is also a lower bound, at least when $n-d(x,y)$
is large enough. As a consequence, we get the same upper and lower
bound for the Green function, away from the diagonal. In rank one,
we recover in a simpler way the main result of Lalley \cite{L1},
specialized to the simple random walk (see \cite{AST} for more
details).

Our method consists in analyzing carefully the transition density,
using the inverse Fourier transform. Recall that Fourier analysis
was developed in the seventies by Macdonald \cite{M1} for $p$--adic
like buildings. It was resumed recently, first by Cartwright
\cite{C} for affine buildings of type $\widetilde{A}_r$ and next
by Parkinson (\cite{P1}, \cite{P2}) in the general case. Notice
that these authors used it already to study isotropic random walks
(\cite{CW}, \cite{P1}, \cite{P3}). They obtained in particular
local and central limit theorem i.e. the asymptotics of the
transition densities $p^n(x,y)$ when $n\to+\infty$  and $x,y$
remain fixed.

Our paper is organized as follows. In Section 2, we recall the
setting of our study and specify the basic objects involved:
affine buildings (of type $\widetilde{A}_r$), the (inverse)
Fourier transform and the ``simple'' random walk on these spaces.
Section 3 is devoted to the rank 2 case, which is typical of the
higher rank case and which is easier to deal with first. In this
case, our method works in fact for any isotropic nearest neighbor
random walk. Moreover we obtain the same upper and lower bound for
$p^n(x,0)$ in the full range $|x|\le n$. Section 4 deals with the
general case. The result is similar, except that the lower bound
is not shown to hold in the range $n-C\le|x|\le n$, where $C$ is
some positive constant (possibly large). In Section 5, we deduce
sharp estimates (same upper and lower bound) for the Green
function, at or above the bottom of the $l^2$ spectrum.

\section{Preliminaries}
\subsection{Root system}
Let $R$ be a root system of type $A_r$ in a real vector space
$\ka$. Let $\ka_\C=\ka+i\ka$ be the complexification of $\ka$,
equipped with its inner product $\bra \cdot, \cdot \ket$. We shall
briefly introduce some standard notation (for more details see e.g
\cite{B}). First let $R^+$ be a choice of positive roots. Let
$\{\alpha_1,\dots,\alpha_r\}$ be the set of simple roots. We
denote by $\ka_+$ the associated positive Weyl chamber and by
$\overline{\ka_+}$ its closure. Let
$\{\lambda_1,\dots,\lambda_r\}$ be the set of fundamental weights.
Let $P=\sum_{i=1}^r \Z\la_i$ be the lattice of weights. Let $P^+$
be the subset of dominant weights, i.e. which lie in
$\overline{\ka_+}$ and let $P^{++}$ be the subset of strictly
dominant weights i.e., those which lie in $\ka_+$. Let
$Q=\sum_{i=1}^r\Z\alpha_i$ be the lattice of roots. The lattice $P$ is the set of vertices of a simplicial complex,
which is called the Coxeter complex. We denote by
$W_0$ the Weyl group and by $\widetilde{W}$ the extended affine
Weyl group (see e.g \cite{P1}). For $\alpha \in R$, let
$\alpha^\vee=\frac{2}{|\alpha|^2}\alpha$ be the coroot associated
to $\alpha$. We have
$$\frac{1}{2}\sum_{\alpha \in R^+}\alpha^\vee =\sum_{i=1}^r \la_i.$$

Let $q\ge 2$ be an integer. We can define $q_w$ for all $w\in \tilde{W}$. Then if $t_\la \in
\widetilde{W}$ is the translation by $\la$, we have (see e.g.
\cite{P2})
$$
q_{t_\la}=q^{\sum_{\alpha \in R^+}\bra\la,\alpha\ket}.
$$
If $\la\in P^+$, we denote by $W_{0\la}$ the stabilizer of $\la$
under the action of $W_0$. If $\la=\sum_{i=1}^r n_i\la_i$, with
$n_i\in \N$ for all $i$, we denote by $|\la|=\sum_{i=1}^r n_i$ the
length of $\la$.
Eventually, the function $\pi$ is defined on $P$ by
$$
\pi(\la)=\prod_{\alpha \in R^+}\bra\alpha^\vee,\la\ket.
$$
\subsection{The symmetric Macdonald polynomials}
The Weyl denominator $\Delta$ and the functions $\mathbf{c}$ and $b$ are defined respectively, for $z \in \ka_\C$, by
\begin{gather*}
\Delta(z)=\prod_{\alpha \in R^+}(e^{\frac{\bra\alpha^\vee,z\ket}{2}}-e^{-\frac{\bra\alpha^\vee,z\ket}{2}}),\\
\mathbf{c}(z)=\prod_{\alpha\in R^+}\frac{1-q^{-1}e^{-\bra\alpha^\vee,z\ket}}{1-e^{-\bra\alpha^\vee,z\ket}},\\
\frac{1}{\mathbf{c}(z)}=\Delta(z)b(z)e^{-\sum_{i=1}^r\bra\la_i,z\ket}.
\end{gather*}
In particular, $|b(i\teta+s)|$ is bounded above and below by a fixed strictly positive constant, for
$(\teta,s)\in i\ka\times \overline{\ka_+}$. The symmetric Macdonald polynomial is defined (see \cite{M2}, \cite{P2})
for $\la \in P^+$ and $z \in \ka_\C$ by
\begin{equation}
\label{sphericalpol}
P_\la(z)=\frac{q_{t_\la}^{-{\frac{1}{2}}}}{W_0(q^{-1})}\sum_{w \in W_0}\mathbf{c}(w \cdot z) e^{\bra \la, w \cdot z \ket}.
\end{equation}
where $\cdot$ denotes the action of $W_0$ on $\ka_\C$. Moreover, for the root systems of type $A_r$, we have
\begin{equation}
\label{sphericalpol2}
P_{\la_i}(z)=q_{t_{\la_i}}^{\frac{1}{2}}\frac{1}{N_{\lambda_i}}\sum_{\lambda \in W_0 \cdot \la_i}e^{\bra \la, z \ket}.
\end{equation}
We define the function $h$ for $z\in \ka_\C$ by
$$
h(z)=\sum_{i=1}^r\sum_{\lambda \in W_0 \cdot \la_i}e^{\bra \la , z \ket}.
$$
\subsection{Affine building and averaging operators}
An affine building (see \cite{R} or \cite{P1}) of type $\tilde{A}_r$ is a nonempty simplicial complex containing subcomplexes
called apartments such that:
\begin{itemize}
\item[$\bullet$] Each apartment is isomorphic (see \cite{P1}) to the Coxeter complex.
\item[$\bullet$] Given two chambers (simplices of maximal dimension) there is an apartment containing both.
\item[$\bullet$] Given two apartments that contain at least a common chamber, there exists a unique isomorphism between them,
which fixes pointwise their intersection.
\end{itemize}
The building will be assumed to be regular. By definition this
means that given any chamber $C$ and any face $F$ (simplex of
codimension $1$) of $C$, the cardinality of the set of chambers
different from $C$ and containing $F$ is independent of $C$ and
$F$, and is equal to $q$. We denote by $\kX$ the set of vertices
(simplices of dimension $1$) of the building. Observe for
instance, that if $r=1$, the building is a tree such that each
vertex has $q+1$ neighbors.

\vspace{0.3cm} We fix a vertex $0$ called the origin, an apartment
$A_0$ containing $0$, and a chamber $C_0$ of $A_0$ containing $0$.
We can identify the set of vertices of $A_0$ with the elements of
the weight lattice $P$. Then we can identify $P^+$ with a subset
of $A_0$ containing the vertices of $C_0$. This subset is denoted
by $A_0^+$. Now given $x \in \kX$, there exists an apartment
containing $C_0$ and $x$. There is also an isomorphism between
this apartment and $A_0$ fixing $C_0$. The image of $x$ by this
isomorphism has a unique conjugate by the Weyl group $W_0$ which
lies in $A_0^+$. This conjugate is called the radial part or
coordinate of $x$ and is identified with an element of $P^+$.

For $\la \in P^+$, we denote by $V_\la(0)$ the sphere of radius
$\la$ around $0$, which by definition is the set of vertices in
the building with radial part equal to $\la$. Its cardinality is
denoted by $N_\la=|V_\la(O)|$. Then we have (see \cite{C})
$$
N_\la=\frac{W_0(q^{-1})}{W_{0\la}(q^{-1})}q_{t_\la},
$$
where $V(q^{-1})=\sum_{w\in V}q_w^{-1}$, for all subgroups $V$ of
$W_0$. For $x\in V_\la(O)$ we set $|x|=|\la|$ and
$\overline{x}=\la$. We set also $x_i=\bra \alpha_i,\la\ket$ for
all $i\le r$.

\vspace{0.5cm}
By $\mathcal{A}$ we denote the
algebra of averaging symmetric operators on $\kX$. It was proved
in \cite{C} that $\mathcal{A}$ is a commutative algebra generated
by the operators
$$
\triangle_j f(x)=\frac{1}{N_{\la_j}} \sum_{y \in V_{\la_j}(O)}
f(y), \quad j=1,\dots,r,
$$
where $f$ is a complex-valued function on $\kX$. Consider $l^2(\kX)$ with a natural scalar product
$$
\bra f, g \ket = \sum_{x \in \kX} f(x) \overline{g(x)}.
$$
Then the closure $\overline{\mathcal{A}}$ is a commutative
$C^*$-algebra. Moreover, $\overline{\mathcal{A}}$ is isometrically
isomorphic to the algebra of $W_0$-invariant continuous functions
on
$$
U=\{\theta \in \ka \mid \text{for all } \alpha \in R,\
\bra\alpha,\teta\ket\ \le \pi\},
$$
and the Gelfand map is given by
$$
\widehat{\triangle_j} = P_{\la_j}.
$$
We observe here that $U$ is $W_0$-invariant and a fundamental domain for the action of the lattice $2\pi Q$ on $\ka$.
Eventually, for $A \in \overline{\mathcal{A}}$ we have the inversion formula
$$
A \delta_y(x)=\frac{W_0(q^{-1})}{|W_0|}
              \int_U \widehat{A}(\theta) \overline{\widehat{\triangle_j}(\theta)} \frac{d\theta}{|\mathbf{c}(\theta)|^2},
$$
for $x,y \in \kX$ and $y \in V_{\la_j}(x)$.
\subsection{The simple random walk}
This is defined as the Markov chain on $\kX$, with transition
probabilities given by
$$
p(x,y)=\begin{cases}
        q_{t_{\la_i}}^{-{\frac{1}{2}}}\rho  & \text{ if } y \in V_{\lambda_i}(x) \\
                                     0  & \text{ otherwise,}
       \end{cases}
$$
where
$$
\rho=\frac{1}{\sum_{i=1}^rq_{t_{\la_i}}^{-{\frac{1}{2}}}N_{\la_i}}.
$$
Let also $\tilde{\rho}=\rho h(0)$ be the associated spectral
radius. For example for the tree, i.e. the $\tilde{A}_1$ case, we
have $\tilde{\rho}=\frac{2\sqrt{q}}{q+1}$. In the case
$\tilde{A}_2$ we have
$$
\tilde{\rho}=\frac{3q}{q^2+q+1},
$$
and in the case $\tilde{A}_3$,
$$
\tilde{\rho}=\frac{14 q^2}{(1+q^2)[(q^2+q+1)+2q^{\frac{1}{2}} (q+1)]}.
$$
\subsection{The function $F_0$}
It is defined on $P^+$ by
$$
F_0(\la)=P_\la(0).
$$
The following Proposition is the analogue of a result obtained in \cite{A} and generalized in \cite{S}.
\begin{prop}
In $P^+$,
\begin{equation}
\label{eq1_5_30}
F_0(\la) \asymp\footnote{we say that $f\asymp g$, when there exists a constant $C>0$ such that, $\frac{1}{C}g(\la)\le f(\la) \le Cg(\la)$
for all $\la$.} q_{t_\la}^{-\frac{1}{2}}\prod_{\alpha \in R^+}(1+\bra\alpha^\vee,\la\ket).
\end{equation}
Moreover,
\begin{equation}
\label{eq2_5_30}
F_0(\la)\sim\footnote{we say that $f\sim g$, when $\frac{f}{g}\to 1$.} \text{const} \cdot \pi(\la)q_{t_\la}^{-\frac{1}{2}}.
\end{equation}
when $\bra\alpha,\la\ket \to +\infty$, for all $\alpha \in R^+$.
\end{prop}
\begin{proof} First we multiply \eqref{sphericalpol} by
$\pi(i\teta)$ removing the singularities of the
$\mathbf{c}$-function, and then we apply the operator
$\left.\pi(-i \partial)\right|_{\teta=0}$. Since the left hand
side is equal to $F_0(\lambda)$ up to a constant, we get
$$
q_{t_\lambda}^{\frac{1}{2}} F_0(\lambda) = p(\lambda),
$$
where $p$ is a polynomial in coordinates of $\lambda$ with highest
order term proportional to $\pi(\la)$. This proves
\eqref{eq2_5_30} and \eqref{eq1_5_30} away from the walls. Next we
extend our estimate along the walls by using a local Harnack
principle. This is obtained immediately by using that $F_0$ is an
eigenfunction of the averaging operators:
$$
\sum_{y \in V_{\la_i}(x)}p(x,y)F_0(\overline{y})=\rho |W_0\cdot \la_i| F_0(\overline{x}),
$$
for all $x\in \kX$ and all $i\le r$.
\end{proof}

\section{Heat kernel estimates: the case $\tilde{A}_2$}
Let $n\in \N$ and $x\in V_\la(O)$. Let
$\alpha_0=\alpha_1+\alpha_2$. We set $\delta=\frac{1}{n+2}
(\la+\la_1+\la_2)$. For $i=0, 1, 2$ we put $\delta_i=\bra \delta,
\alpha_i \ket$. Let
$$
\phi(\delta) = \min\{u \in \overline{\ka_+}\mid \log h(u) - \bra\delta, u\ket \}.
$$
The main goal of this section is to prove
\begin{theo}
\label{theoestimateheatA2}
The following estimate
\begin{equation}
\label{estimateheatA2}
p^n(0,x) \asymp \frac{1}{n^3}\rho^n e^{n\phi(\delta)} F_0(\overline{x})
                  \frac{1}{\sqrt{n^2(1-\delta_0)(1-\delta_1)(1-\delta_2)}},
\end{equation}
holds uniformly on the set $\{|x| \le n-1\}$.
\end{theo}
We will see that the function $e^{\phi(\delta)}$ is bounded. Thus the exponent $n$ in the theorem can be replaced by $n+2$, which appears more naturally
in the proof.
The next theorem gives a more precise statement of the estimate at the boundary of the domain.
We adopt the following notation for the binomial coefficients:
$$
C^n_k = \frac{n!}{k!(n-k)!}.
$$
\begin{theo}
\label{theo2estimateheatA2}
Let $K>0$. Then
\begin{equation}
\label{asymptoticboundary}
p^n(0,x)\asymp n^d (\rho q^{-1})^n C^{n-d}_{x_1\vee x_2-d}
\end{equation}
uniformly in the set $\{n\ge |x| \ge n-K\}$, where $d=n-|x|$.
\end{theo}
\begin{note}
When $n-x_1\vee x_2\le K'$ for some fixed constant $K'>0$, then the estimate becomes
$$
p^n(0,x)\asymp (\rho q^{-1})^n n^{(n-x_1\vee x_2)+d}.
$$
\end{note}

\begin{rem}
\emph{Here is the corresponding result for the tree (cf
\cite{AST}). In this case we can give an explicit formula of the
function $\phi$ appearing in the estimate. In fact we have
$$p^n(0,x)\asymp
\frac{|x|}{n\sqrt{n-|x|}}\rho^ne^{n\phi(\delta)}q^{-\frac{|x|}{2}},$$
where
$$\phi(\delta)=\frac{1}{2}\{(1+\delta)\log(1+\delta)+(1-\delta)\log(1-\delta)\}.$$
}
\end{rem}

\subsection{Proof: the beginning}
If $\teta=\teta_1\alpha_1+\teta_2\alpha_2$, we set
$|\teta|_\infty=\max\{|\teta_1|,|\teta_2|\}$. We say that a weight
is away from a wall, when its distance to the wall is larger than
some fixed constant (which will be determined in the proof). We
denote by $C$ a constant whose value may change from line to line.

We begin by some elementary transformations of $p^n(O,x)$. First, by using \eqref{sphericalpol2}, we get
\begin{eqnarray*}
&p^n(0,x)&=C\int_U \left(\frac{1}{2}P_{\la_1}(i\teta)+\frac{1}{2}P_{\la_2}(i\teta)\right)^nP_\la(i\theta)
                                                  \frac{d\theta}{|\mathbf{c}(i\theta)|^2}\\
        &&=C\rho^n\int_U h^n(i\theta)P_\la(i\theta)\frac{d\theta}{|\mathbf{c}(i\theta)|^2}.
\end{eqnarray*}
Next by \eqref{sphericalpol} and the $W_0$-invariance of $h$, we
get
\begin{equation}
\label{heatformula}
p^n(0,x)=C\rho^nq_{t_\la}^{-{\frac{1}{2}}}\int_U h^n(i\theta)e^{-i\bra\teta,\la+\la_1+\la_2\ket}\Delta(i\teta)b(i\teta)d\teta.
\end{equation}
Now we make two elementary observations. First
\begin{equation}
\label{productformula}
h+2=(1+e^{\la_1})(1+e^{-\la_2})(1+e^{\la_2-\la_1}).
\end{equation}
Moreover, we have the following lemma, whose proof is left to the
reader (see also Section \ref{sectionAr} for a more general
result).
\begin{lem}
\label{lemderiv}
$$
\pi(\partial)[h^{n+3}]=(n+3)(n+2)(n+1)\left[\frac{n+3}{n+1}h+2\right]h^n\Delta.
$$
\end{lem}
The idea now is to use Lemma \ref{lemderiv} and integrate by parts
in \eqref{heatformula}. This will be done in two different ways,
depending on if $|x|\le \frac{n}{2}$ or if $|x|> \frac{n}{2}$. We
notice that the factor $\frac{1}{2}$ plays no role here, and could
be replaced by $1-\eta$ for any $\eta \in (0,1)$.

\subsection{The case when $|x|\le \frac{n}{2}$}
We consider the function
$$
q^n(x)=\frac{n+3}{(n+1)\rho}p^{n+1}(0,x)+2p^n(0,x).
$$
By Lemma \ref{lemderiv}, 
after an integration by parts we get the expression
$$
q^n(x)=C\frac{\rho^nq_{t_\la}^{-{\frac{1}{2}}}}{(n+3)(n+2)(n+1)}\int_U h^{n+3}(i\theta)\pi(\partial)
    \left[e^{-i\bra\teta,\la+\la_1+\la_2\ket}b(i\teta)\right]d\teta.
$$
Now we make the change of variables $i\teta\to i\teta+s$ for some $s\in \overline{\ka_+}$ whose value will be specified in
the sequel, and we find
$$
q^n(x)=C\frac{\rho^nq_{t_\la}^{-{\frac{1}{2}}}e^{-\bra \la+\la_1+\la_2,s \ket}}{(n+3)(n+2)(n+1)}
        \int_U h^{n+3}(i\theta+s)\pi(\partial)[e^{-i\bra\teta,\la+\la_1+\la_2\ket}b(i\teta+s)]d\teta.
$$
Moreover $b(i\teta+s)$ and all its derivatives are bounded functions of $(\teta,s)\in U\times\overline{\ka_+}$. Thus for $x$ or
$\la$ sufficiently away from the walls,
\begin{eqnarray*}
q^n(x)&=C\frac{\rho^nq_{t_\la}^{-{\frac{1}{2}}}e^{-\bra\la+\la_1+\la_2,s\ket}\pi(\la+\la_1+\la_2)}{(n+3)(n+2)(n+1)}\\
      &\times\int_U h^{n+3}(i\theta+s)e^{-i\bra\teta,\la+\la_1+\la_2\ket}b(i\teta+s)d\teta + \text{lower order terms}.
\end{eqnarray*}

\subsection{The case when $|x|>\frac{n}{2}$}
\label{range2}
In this case $q^n$ and $p^n$ are not anymore comparable. However, when $\la$ is sufficiently away from the wall
$\{\alpha_1-\alpha_2=0\}$, then it results from \eqref{productformula} and our choice of $s$ in the next
subsection (see also Remark \ref{remshift}), that at least for $n$ sufficiently large, the function
$\teta \mapsto \frac{1}{\frac{n+3}{n+1}h(i\teta+s)+2}$ does not vanish on $U$. Thus after an integration by parts in \eqref{heatformula},
we see that for $\la$ away from the walls $\{\alpha_2=0\}$ and $\{\alpha_1-\alpha_2=0\}$,
\begin{eqnarray*}
p^n(0,x)&=&C\frac{\rho^nq_{t_\la}^{-{\frac{1}{2}}}e^{-\bra\la+\la_1+\la_2,s\ket}\pi(\la+\la_1+\la_2)}{(n+3)(n+2)(n+1)}\\
        &\times& \int_U h^{n+3}(i\theta+s)e^{-i\bra\teta,\la+\la_1+\la_2\ket}(b_1+b_2)(i\teta+s)d\teta,
\end{eqnarray*}
where $b_1=\frac{b}{\frac{n+3}{n+1}h+2}$ and
$b_2=\frac{e^{i(\la+\la_1+\la_2)}}{\pi(\la+\la_1+\la_2)}\pi(\partial)[b_1e^{-i(\la+\la_1+\la_2)}]-b_1$
is the remainder.

\subsection{Choice of $s$ and the stationary phase method}
In the preceding subsections we have seen that $q^n(x)$, in the range $|x|\le \frac{n}{2}$, and $p^n(0,x)$, in the range
$|x|>\frac{n}{2}$, were comparable for $\la$ away from the walls to
\begin{equation}
\label{maintermheat}
C(n,\la)\int_U \frac{h^{n+2}(i\theta+s)}{h^{n+2}(s)}e^{-i\bra \teta,\la+\la_1+\la_2\ket}\tilde{b}(i\teta+s)d\teta,
\end{equation}
with
$$
C(n,\la)=\frac{\rho^nq_{t_\la}^{-{\frac{1}{2}}}e^{-\bra \la+\la_1+\la_2,s \ket}h^{n+2}(s)\pi(\la)}{n^3},
$$
and $\tilde{b}$ equal to $hb$ or $h(b_1+b_2)$ respectively if $|x|\le \frac{n}{2}$ or if $|x|>\frac{n}{2}$. Now we choose the
shift $s$ according to the stationary phase method, i.e. as a solution in $\overline{\ka_+}$ of the equation
\begin{equation}
\label{equationshift}
\frac{\nabla h}{h}(s)=\delta.
\end{equation}
Solving such an equation is classical. We consider the function
$\phi^\delta :u\mapsto \log h(u)-\bra\delta,u\ket$ on $\ka$. Since
$|\delta|<1$ (by hypothesis $|x|<n$), this function tends to
infinity when $|u|\to +\infty$. Moreover, $h$ is $W_0$-invariant,
and $\bra \delta,u \ket$ is maximal when $u\in \overline{\ka_+}$.
Thus $\phi^\delta$ attains its minimum in $\overline{\ka_+}$ at
some point $s=s(\delta)$, which satisfies equation
\eqref{equationshift}. Without loss of generality, we will assume
in the sequel that $\bra \la_1,s\ket \ge \bra\la_2,s\ket$. The
following lemma collects some properties of $s$.
\begin{lem}
\label{propertiesshift}
\begin{enumerate}
\item The condition $\bra\lambda_1,s\ket \ge \bra\lambda_2,s\ket$ implies $\delta_1\ge \delta_2$.
\item The function $\phi^\delta$ is strictly convex, thus $s$ is the unique point where it attains its minimum,
and $\delta\mapsto s$ is continuous on $\{|\delta|<1\}$.
\item When $|\delta| \to 1$, $|s|\to +\infty$, where $|s|$ is the Euclidean norm of $s$. More precisely, when $|\delta|\to 1$,
$$
e^{\bra \la_2-\la_1,s \ket} = \frac{1-\delta_1}{\delta_1}+O(e^{-\bra \la_2,s \ket}),
$$
$$
e^{-\bra \la_2,s \ket}\asymp 1-|\delta|.
$$
\item We have $s=0$ if, and only if, $\delta=0$. When $\delta_1-\delta_2\to 0$, then
$$
\bra\la_2-\la_1,s\ket \asymp \delta_1-\delta_2.
$$
\end{enumerate}
\end{lem}
\begin{proof} Equation \eqref{equationshift} is equivalent to
the system
\begin{eqnarray}
\label{systemequationshift}
\left\{\begin{array}{ccc}
             \sinh \bra\la_1,s\ket + \sinh \bra\la_1-\la_2,s\ket =& \delta_1 h\\
             \sinh \bra\la_2,s\ket - \sinh \bra\la_1-\la_2,s\ket =& \delta_2 h.
          \end{array}\right.
\end{eqnarray}
The first statement follows since $\sinh$ is increasing and has
the same sign than its argument. The convexity of $\phi^\delta$
comes from a more general result: assume that $h=\sum_\la e^\la$
is a sum of exponentials. Then the second order derivative
$d^2\phi^\delta_u$ of $\phi^\delta$ at some point $u\in \ka$ is
given by
$$
d^2\phi_u^\delta= \frac{hd^2 h-|\nabla h|^2}{h^2}(u)=
\frac{1}{h^2(u)}\sum_{\la\neq \la'}[\la-\la']^2e^{\bra\la+\la',u\ket},
$$
which implies that $d^2\phi_u^\delta$ is positive definite. Thus
$\phi^\delta$ is strictly convex. The second point follows
immediately. By adding the both equations of
\eqref{systemequationshift}, and because $\sinh < \cosh$, we see
that $|s|\to \infty$ when $|\delta|\to 1$. Multiplying by
$e^{-\bra \la_1,s\ket}$ in both side of the first equation gives
immediately the asymptotic of $e^{\bra \la_2-\la_1,s\ket}$. Doing
the same in the sum of the two equations gives the asymptotic of
$e^{-\bra\la_2,s\ket}$. By doing now the difference between them,
we obtain in the same way the last point. The assertion that $s=0$
if, and only if, $\delta =0$ is straightforward. This concludes
the proof of the lemma.
\end{proof}

\begin{rem}
\label{remshift}
\emph{\begin{enumerate}
\item As announced, the estimates of the lemma show that the function $e^{\phi(\delta)}=h(s)e^{-\bra\delta,s\ket}$ in
      \eqref{estimateheatA2} is bounded.
\item The last point of the lemma implies that $(1+e^{\bra\la_2-\la_1,i\teta+s\ket})$ is larger (up to a constant) than
$\delta_1-\delta_2$ for any $\teta \in U$. Thus, thanks to Formula
(\ref{productformula}), we see that for $n$ and
$\bra\alpha_1-\alpha_2,\la\ket$ large enough, $\frac{n+3}{n+1}h+2$
does not vanish in $U$. This justifies our assumption of section
\ref{range2}.
\end{enumerate}}
\end{rem}
We consider now the phase function
\begin{eqnarray}
\label{phase} F^\delta(\teta)=\log h(i\teta +s)-\log
h(s)-i\bra\delta,\teta\ket
\end{eqnarray}
which is well defined at least in a small neighborhood of $0$,
independent of $s$ (or $\delta$). By our choice of $s$,
$F^\delta(0)=0$ and $\nabla F^\delta(0)=0$. The next lemma gives a
much more precise result on the behavior of $F^\delta$ near $0$.
Let $\Re F^\delta$ and $\Im F^\delta$ denote respectively the real
and imaginary part of $F^\delta$.
\begin{lem}
\label{quadraticformphase}
There exists two constants $\epsilon>0$ and $C>0$ such that,
$$
\Re F^\delta(\teta) \asymp -q_\delta(\teta),\ \text{and}\ |\Im
F^\delta(\teta)|\le C |\teta|_\infty q_\delta(\teta),
$$
uniformly for $|\teta|_\infty \le \epsilon$ and $|\delta| < 1$, where
$$
q_\delta(\teta)=e^{\bra\la_2-\la_1,s\ket}\bra\la_1-\la_2,\teta \ket^2+e^{-\bra\la_2,s\ket}\bra\la_2,\teta\ket^2.
$$
\end{lem}
\begin{proof} Since $F^\delta(0)=0$ and $\nabla F^\delta(0)=0$,
we see that for all $\teta \in U$, there exists $\teta' \in U$
such that $|\teta'|_\infty \le |\teta|_\infty$ and
$$
F^\delta(\teta)=d^2F^\delta_{\teta'}(\teta),
$$
where $d^2F^\delta_{\teta'}$ is the second order derivative of $F^\delta$ in
$\teta'$. We can compute it, as we did for $\phi^\delta$ in Lemma
\ref{propertiesshift}:
\begin{eqnarray*}
d^2F^\delta_{\teta'}(\teta)&=&\frac{-1}{h(i\teta'+s)^2}\sum_{\la,\la'}\bra\la-\la',\teta\ket^2e^{\bra\la+\la',i\teta'+s\ket}\\
                   &=&\frac{-1}{|h(i\teta'+s)|^4}\sum_{\la,\la',\mu,\mu'}\bra\la-\la',\teta\ket^2e^{\bra\la+\la'+\mu+\mu',s\ket}
               e^{i\bra\la+\la'-\mu-\mu',\teta'\ket}.
\end{eqnarray*}
Next observe that $|h(i\teta'+s)e^{-\bra\la_1,s\ket}|$ is bounded above and below by strictly positive constants,
when $|\teta'|_\infty$ is small. So if we multiply up and down the right member of the last equality by
$e^{-4\bra\la_1,s\ket}$, and then take successively the real part and the imaginary part, we get the two assertions of
the lemma.
\end{proof}
We denote by $J(n,\la)$ the integral appearing in \eqref{maintermheat}. According to the notation of the preceding
lemma, let $\epsilon'<\epsilon$ (taken small in the next proposition). We divide $J(n,\la)$ into the sum of the integral,
let say $J_1(n,\la)$, over $[-\epsilon',\epsilon']^2$ and the integral $J_2(n,\la)$ over
$U\smallsetminus [-\epsilon',\epsilon']^2$, where $[-\epsilon',\epsilon']^2$ denotes the set $\{|\teta|_\infty \le \epsilon'\}$.
\begin{prop}
There exists $\epsilon'>0$ (independent of $\delta$ and $n$), such that the following estimates hold for $n$ large enough
and $\la$ away from the walls
$$
|J_1(n,\la)| \asymp \frac{1}{\sqrt{n^2(1-\delta_1)(1-|\delta|)}},
$$
$$
|J_2(n,\la)| \le \frac{C}{\sqrt{n(1-\delta_1)}}e^{-cn(1-|\delta|)},
$$
where $C$ and $c$ are two strictly positive constants.
\end{prop}
\begin{proof} We have
$$
J_1(n,\la)=\int_{|\teta|_\infty \le
\epsilon'}e^{(n+2)F^\delta(i\teta+s)}\tilde{b}(i\teta+s)d\teta.
$$
Thus essentially the first statement of the proposition is given by
Lemma \ref{quadraticformphase}, and a change of variable. In fact
we just need in addition a control of $\tilde{b}(i\teta+s)$ for
small $\teta$'s. In the case where $|x|\le \frac{n}{2}$, this is
immediate, since $|b|$ is bounded above and below by strictly
positive constants, and by continuity of $\delta \mapsto s$, $|h|$
also, if $\epsilon'$ is sufficiently small. In the other case,
where $|x|>\frac{n}{2}$, by \eqref{productformula}, we see that
$|h+2|$ stays away from $0$ if $|\teta|_\infty$ is small. Thus for
$\la$ sufficiently away from the walls, $b_2$ becomes negligible
in front of $b_1$. Then since $(hb_1)(s)$ is real and bounded
(above and below) on $\overline{\ka_+}$, this gives the desired
estimate of $|J_1(n,\la)|$. The estimate of $J_2$ is more
complicated, since $|hb_2|$ may become larger than $|hb_1|$ and
even explode, for instance when $\bra\la_2-\la_1,\teta\ket=\pm
\pi$, and $\delta_1-\delta_2$ tends to $0$ (cf.
\eqref{productformula}). Fortunately, as we will see, this is
compensated by the exponential decay of $e^{(n+2)\Re
F^\delta(i\teta+s)}$. Indeed
$$
\frac{|h(i\teta+s)|}{h(s)}=e^{\frac{1}{2}\log(1-(1-\frac{|h(i\teta+s)|^2}{h^2(s)}))}\le
                           e^{-\frac{1}{2}[1-\frac{|h(i\teta+s)|^2}{h^2(s)}]}.
$$
But
$$
1-\frac{|h(i\teta+s)|^2}{h^2(s)}=\frac{1}{h^2(s)}\sum_{\la,\la'}e^{\bra\la+\la',s\ket}(1-\cos\bra\la-\la',\teta\ket).
$$
Multiplying again the numerator and denominator by
$e^{-2\bra\la_1,s\ket}$, we obtain
$$
1-\frac{|h(i\teta+s)|^2}{h^2(s)}\ge c\sum_{\la,\la'}e^{\bra\la+\la'-2\la_1,s\ket}\bra\la-\la',\teta\ket^2,
$$
for some constant $c>0$. Observe now that since $s=0 \Leftrightarrow \delta=0$ and $s$ is continuous, then $s$ stays away from
$0$ when $|x|>\frac{n}{2}$. In particular $\bra\la_1,s\ket$ and $\bra\la_2,s\ket$ stay also away from $0$. Therefore with
\eqref{productformula} we see that $|b_1|$ or $|b_2|$ can explode only when $1+e^{\la_2-\la_1}$ is small, i.e when
$\bra\la_2-\la_1,\teta\ket\pm \pi$ and $\delta_1-\delta_2$ are small. But in this case
$$
(n+2)e^{\bra\la_2-\la_1,s\ket}\bra\la_2-\la_1,\teta\ket^2\asymp (n+2)(1-\delta_1)\asymp n,
$$
since $\delta_1-\delta_2$ small implies that $1-\delta_1$ is away from $0$. As a consequence
$$
e^{-\frac{c}{2}(n+2)\sum_{\la,\la'}e^{\bra\la+\la'-2\la_1,s\ket}\bra\la-\la',\teta\ket^2}(|hb_1(i\teta+s)|+|hb_2(i\teta+s)|)
$$
is bounded in $U$. The estimate of $|J_2(n,\la)|$  follows with Lemma \ref{propertiesshift} and a change of variable.
\end{proof}
By the preceding proposition, when $n(1-|\delta|)\to \infty$,
$J_2$ becomes negligible in front of $J_1$. Thus we can find a
constant $K>0$ such that $C(n,\la)|J(n,\la)|$ has the estimate
\eqref{estimateheatA2} when $\la$ is away from the walls and
$n(1-|\delta|)\ge K$. Moreover the preceding proposition implies
also that $C(n,\la)|J(n,\la)|$ is bounded by the expression in
\eqref{estimateheatA2}. Now the rest of the proof can be
decomposed into two steps. First we prove the lower estimate when
$n(1-|\delta|)\le K$. Then we extend our estimate along the walls
by using a local Harnack principle, and we prove by the way that
$q^n(x)$ is comparable to $p^n(0,x)$ when $|x|\le \frac{n}{2}$.
%
%

\subsection{Lower bound when $n(1-|\delta|)\le K$}
\label{lowerboundA2boundary} We present two proofs. The first is analytical, and the second is purely combinatorial.
One advantage of the second proof is that it is valid in the whole range $n(1-|\delta|)\le K$, whereas the first is only
valid when $\la$ is away from the walls $\{\alpha_2=0\}$ and $\{\alpha_1-\alpha_2=0\}$. Another advantage of the combinatorial
approach is that it provides a more elementary proof of the upper bound.
\subsubsection{Analytical proof}
We begin by
\begin{lem}
When $(n+2)(1-\delta_1)\to +\infty$ and $(n+2)(\delta_1-\delta_2) \to +\infty$, then $\tilde{b}(i\teta+s)\to 1$,
uniformly in $U$.
\end{lem}
\begin{proof} Two cases may cause problems. Either when
$\delta_1-\delta_2$ tends to $0$. But in this case
$\bra\alpha,\la\ket\asymp n$ for all $\alpha \in \R^+$, since the
line $\la_1=\la_2$ does not cross walls in the range
$n(1-|\delta|) \le K$ (at least if $n$ is large enough). Thus when
$n(\delta_1-\delta_2)\to +\infty$, $|b_2|$ tends to $0$. The other
case is when $(1-\delta_1)\to 0$, because we are not sure a priori
that the function $b$ tends to $1$. But this becomes true if
$n(1-\delta_1) \to +\infty$. Now the proof of the lemma is
straightforward.
\end{proof}
By the preceding lemma there exists a constant $K'>0$ such that if $(n+2)(1-\delta_1)\ge K'$ and
$(n+2)(\delta_1-\delta_2)\ge K'$, then $J(n,\la)$ is comparable to
$$
\frac{1}{h^{n+2}(s)}\int_U h^{n+2}(i\teta+s)e^{-i\bra\la+\la_1+\la_2,\teta\ket}d\teta.
$$
We denote by $I(n,\la)$ the integral in the last expression. We compute it by developing $h^{n+2}(i\teta+s)$ and using that
the integral of $e^{i\bra\mu,\teta\ket}$ is null whenever $\mu$ is a non zero weight. We get, if $|U|$ denotes
the Lebesgue measure of $U$,
$$
I(n,\la)=|U|\sum\frac{(n+2)!}{n_1!n_2!\dots n_6!},
$$
where the sum is over the family of integers $(n_1,\dots,n_6)$ such that $n_1+\dots+n_6=n+2$, $n_1-n_2+n_5-n_6=x_1+1$ and
$n_3-n_4-n_5+n_6=x_2+1$. In particular we can choose $n_1=x_1+1-d$, $n_3=x_2+1+d$, $n_5=d$, $n_2=n_4=n_6=0$, with
$d=n-x_1-x_2$. Thus, by Stirling's formula, we get
\begin{eqnarray*}
I(n,\la)&\ge& c_1 \frac{(n+2)^{n+2}}{x_1^{x_1+1-d}x_2^{x_2+1+d}}\frac{1}{\sqrt{x_2}}\\
        &\ge& c_2 n^d (\frac{n+2}{x_1})^{x_1+1-d}(\frac{n+2}{x_2})^{x_2+1+d}\frac{1}{\sqrt{n(1-\delta_1)}}\\
        &\ge& c_3 n^d \frac{1}{\delta_1^{x_1+1-d}}\frac{1}{(1-\delta_1)^{x_2+1+d}}\frac{1}{\sqrt{n(1-\delta_1)}},
\end{eqnarray*}
where $c_1,c_2$ and $c_3$ are strictly positive constants. The passage from the second to the third line is justified by the
inequality $\frac{(n+2)(1-\delta_1)}{x_2}\ge 1$. On the other side, with Lemma \ref{propertiesshift} we get
\begin{eqnarray*}
h(s)^{n+2}e^{-\bra\la+\la_1+\la_2,s\ket} \asymp (\frac{\delta_1}{1-\delta_1})^{x_2+1+d}(1+\frac{1-\delta_1}{\delta_1})^{n+2}
                                          \frac{1}{(1-|\delta|)^{(n+2)(1-|\delta|)}}.
\end{eqnarray*}
Then
\begin{eqnarray*}
h(s)^{n+2}e^{-\bra\la+\la_1+\la_2,s\ket} \asymp (\frac{\delta_1}{1-\delta_1})^{x_2+1+d}\frac{1}{\delta_1^{n+2}}n^d.
\end{eqnarray*}
Thus we have proved the lower bound when $(n+2)(1-\delta_1)\ge K'$ (and $\la$ is away from the walls).

\subsubsection{Combinatorial proof}
\label{comblowbound}
We have just seen above that
$$
(h(s)e^{-\bra\delta,s\ket})^n\frac{1}{\sqrt{n(1-\delta_1)}} \asymp C^{n+2-d}_{x_1+1-d}\ n^d,
$$
for $x$ and $n$ such that $(n+2)(1-|\delta|)\le K$. Now observe that in this range
$$
C_{x_1+1-d}^{n+2-d}\asymp \frac{n}{x_2+1}C^{n-d}_{x_1-d},
$$
and $q_{t_\la}\asymp q^{2n}$. Thus we are lead to prove the estimate of Theorem \ref{theo2estimateheatA2}. In fact it is more
convenient to prove the corresponding result for the radial random walk. If $\overline{p}$ denotes its transition kernel,
then by definition $\overline{p}^n(0,\la)=p^n(0,x)N_\la$ for all $n\in \N$, all $\la \in P^+$, and all $x\in V_\la(0)$.
Since $N_\la$ is comparable to $q_{t_\la}$ we have to prove that
$$
\overline{p}^n(0,\la)\ge c (q\rho)^n C^{n-d}_{x_1-d}\ n^d,
$$
for some constant $c>0$. But for any $\la\in P^+$, $q\rho=\overline{p}(\la,\la+\la_1)$, and when $\la\in P^{++}$, we
have also $q\rho=\overline{p}(\la,\la+\la_2)$. Thus it suffices to prove that, if $n+2-|\la|\le K$, then the number of paths
from $0$ to $\la$ in $P^+$, is comparable to $C^{n-d}_{x_1-d}\ n^d$. This is elementary, and can be seen as follows.
Consider the sequence of increments of the radial random walk up to time $n$, $(\epsilon_1,\dots,\epsilon_n)$, with
$\epsilon_i\in \{\pm \la_1,\pm \la_2, \pm (\la_1-\la_2)\}$ for all $i\le n$. Then choose in arbitrary order, $x_1-d$ terms
equal to $\la_1$, $x_2+d$ terms equal to $\la_2$, and $d$ terms equal to $\la_1-\la_2$. The number of ways to do it
is comparable to $C^{n-d}_{x_1-d}\ n^d$, which proves the lower bound. \newline In fact we can also prove the upper bound
by combinatorial arguments. Indeed, there must be at least $x_1-d$ terms equal to $\la_1$ in the sequence, and at most $x_1+d$.
Otherwise the random walk could not reach the point $\la$ at time $n$. Now if $\la_1$ appears $k\in \{x_1-d,\dots,x_1+d\}$
times, then for the same reason there must be also at least $n-d-k$ terms equal to $\la_2$. But when $k$ is fixed, the number
of sequences satisfying these conditions is bounded (up to a constant) by $C^{n-d}_k\ n^{d}$. Moreover, since
$x_1\ge \frac{n-d}{2}$ (and $d\le K$), there exists a constant $c>0$ such that $C^{n-d}_k\le c C_{x_1-d}^{n-d}$, for all
$k\in \{x_1-d,\dots,x_1+d\}$. This proves the desired upper bound.

\subsection{Local Harnack principle and estimate along the walls}
Let us assume that the estimate proved until now is valid when $\overline{x}$ is at distance at least $D$ from the walls
$\{\alpha_2=0\}$ and $\{\alpha_1-\alpha_2=0\}$. From the heat equation satisfied by $p$, we know that there exists a constant
$C>0$, such that
\begin{equation}
\label{locharnackprinc}
 p^n(O,y)\le Cp^{n+1}(O,x),
\end{equation}
for all neighbors $x$ and $y$. Let now $x\in \kX$ be such that $\overline{x}$ is at distance less than $D$ of the wall
$\{\alpha_2=0\}$ for instance. By repeated applications of \eqref{locharnackprinc} we see that there exist vertices $y_1$ and
$y_2$ such that $\overline{y_1}=\overline{x}+D (\la_2 -\la_1)$, $\overline{y_2}=\overline{x}+D\la_2$, and
$$
cp^{n-D}(0,y_1)\le p^n(0,x)\le Cp^{n+D}(0,y_2),
$$
where $c>0$ and $C>0$ are other constants. Therefore, we only need
to show that our upper estimate of $p^{n+D}(0,y_2)$ is comparable
to our lower estimate of $p^{n-D}(0,y_1)$. This yields to prove
that $e^{n[\phi(\delta_{y_2})-\phi(\delta_{y_1})]}$ is bounded,
where $\delta_{y_i}=\frac{\overline{y_i}+\la_1+\la_2}{n+2}$. But
by an elementary calculus, we see that
$$
\nabla \phi(\delta)=-s.
$$
Hence
$$
\bra\nabla \phi(\delta_{y_1}),\delta_{y_2}-\delta_{y_1}\ket\le 0.
$$
This implies well that
$e^{n[\phi(\delta_{y_2})-\phi(\delta_{y_1})]}$ is bounded. This
proves that our estimate extend near the walls. The only missing
part now is to see that $q^n(x)$ and $p^n(0,x)$ are comparable
when $|x|\le \frac{n}{2}$. This results also from \eqref{locharnackprinc}. Indeed it implies that
$$
p^n(0,x)\le q^n(x)\le Cp^{n+3}(0,x),
$$
and our estimates show that $p^{n+3}$ is comparable to $p^n$ when $|x|\le \frac{n}{2}$.

The proof of the theorems \ref{theoestimateheatA2} and \ref{theo2estimateheatA2} is now finished. The extension of
\eqref{asymptoticboundary} when $|x|=n$ is straightforward.

\begin{rem} \emph{For the rank $2$ case, this proof may in fact be
applied for any isotropic nearest neighbor random walk. Let us detail.
Such random walk has transition density given by $p(x,y)=p_i$ if $y\in
V_{\la_i}(x)$, for $i=1,2$, and $p(x,y)=0$ if $x$ and $y$ are not
neighbors. Now in rank $2$, $q_{t_{\la_1}}=q_{t_{\la_2}}$
and $N_{\la_1}=N_{\la_2}$. Thus we have the formula similar to \eqref{heatformula}
$$
p^n(0,x)=C\rho^nq_{t_\la}^{-{\frac{1}{2}}}\int_U h^n(i\theta)e^{-i\bra\teta,\la+\la_1+\la_2\ket}\Delta(i\teta)b(i\teta)d\teta,
$$
but with $\rho^{-1}=q_{t_{\la_1}}^{-\frac{1}{2}}N_{\la_1}$ and
$h=p_1\sum_{\la \in W_0\la_1}e^{\la}+p_2\sum_{\la\in
W_0\la_2}e^\la$. In particular observe that the spectral radius
$\tilde{\rho}:=\rho h(0)$ is the same for all these random walks.
Now remark that for some constants $c,c'>0$, we have
$$
h+c+cc'^3=c(1+c'e^{\la_1})(1+c'e^{\la_2-\la_1})(1+c'e^{-\la_2}).
$$
Namely $c'=p_2/p_1$ and $c=p_1^2/p_2$. So all the preceding proof can be
applied, and Theorem \ref{theoestimateheatA2}
can in fact be deduced for all the isotropic nearest neighbor random walks.}
\end{rem}

\section{Heat kernel estimate for general buildings of type $\tilde{A}_r$}
\label{sectionAr}
Given $x\in V_\la(0)$ and $n\ge 0$, we set
$\delta=\frac{\la+\sum_{i=1}^r\la_i}{n+r}$, and for $i\le r$, let $\delta_i:=\bra\alpha_i,\delta\ket$. We define again $\phi$ by
$$
\phi(\delta)=\min\{u \in \overline{\ka_+}\mid \log h(u)-\bra u ,
\delta \ket\}.
$$
We have the following result
\begin{theo}
\label{theoestimateheatA}
There exists a constant $K>0$, such that the following estimate holds uniformly in the set $\{|x|\le n-K\}$
\begin{equation}
\label{estimateheatA}
p^n(0,x) \asymp \frac{1}{n^{|R^+|}}\rho^n e^{n \phi(\delta)} F_0(\overline{x})
         \frac{1}{\sqrt{n^r\prod_{\alpha\in R^+}
         (1-\bra\alpha,\delta\ket)}}.
\end{equation}
Moreover, the upper estimate holds in the whole domain $\{|x|<n\}$.
\end{theo}
\begin{proof} The proof follows exactly the same lines as in
rank $2$. After elementary computations we get
$$
p^n(0,x)=C\rho^nq_{t_\la}^{-\frac{1}{2}} \int_U h^n(i\teta)e^{-\bra i\teta,\la+\sum_{i=1}^r\la_i\ket}
           \Delta(i\teta)b(i\teta)d\teta.
$$
The next result is the analogue of
Lemma \ref{lemderiv}
\begin{lem}
\label{lemderivAr}
The following formula holds
$$
\pi(\partial)\left[h^{n+|R^+|}\right]=(n+|R^+|)\dots(n+1)r_n(h)h^n\Delta,
$$
where $r_n$ is a polynomial of the form
$$
r_n(h)=\sum_{k=r}^{|R^+|} c_k(h+2)^{k-r}\frac{h^{|R^+|-k}}{(n+1)\dots(n+|R^+|-k)},
$$
with constants $c_k\in \R$.
\end{lem}
\begin{proof} Let us first show that
\begin{equation}
\label{productformulaAr}
h+2=\prod_{\la\in W_0\la_1}(1+e^\la)=\prod_{i=1}^{r+1}(1+e^{\la_i-\la_{i-1}}),
\end{equation}
with the convention $\la_{-1}=\la_{r+1}=0$. 
We prove it by using a well known symbolic description of the fundamental weights and of their conjugates by $W_0$.
We associate to $\la_1$ the symbol $x_1$, to $\la_2$ the symbol $x_1+x_2$,
and so on until $\la_r$ to which we associate the symbol $x_1+\dots+x_r$, with the rule
$x_1+\dots+x_{r+1}=0$. For instance $x_2+\dots +x_{r+1}$ represents the weight $-\la_1$. Then we have a nice description of
the conjugates of a fundamental weight
$$
W_0\la_k=\{x_{i_1}+\dots +x_{i_k}\mid i_i\neq i_2 \dots \neq i_k\},
$$
for all $k\le r$. With this notation it is now elementary to see that
$$
h+2=\prod_{i=1}^r (1+e^{x_i}),
$$
which gives \eqref{productformulaAr}. Next, for any subset $I$ of $R^+$, we define $\pi^I$ on $P^+$ by
$$
\pi^I(\la)=\prod_{\alpha \in I}\bra \alpha^\vee,\la\ket.
$$
We have
$$
\pi(\partial)\left[h^{n+|R^+|}\right]=\sum_{k=1}^{|R^+|}(n+|R^+|)\dots(n+|R^+|-k+1)h^{n+|R^+|-k}f_k,
$$
where
$$
f_k=\sum_{I_1\cup I_2\cup \dots\cup I_k=R^+}
            \left[\pi^{I_1}(\partial)(h+2)\right]\dots \left[\pi^{I_k}(\partial)(h+2)\right],
$$
for all $k\le |R^+|$. The polynomials $f_k$ are skew-invariant.
Thus they are divisible by $\Delta$. Moreover, $f_k$ is of degree
less or equal to $k$. It implies that $f_k=0$ when $k<r$, since
$\Delta$ is of degree $r$. Thus we may assume $k\ge r$. Now we
observe that there is only $r$ roots $\alpha$ in $R^+$ such that
$\bra \alpha,\la_1\ket\neq 0$. The same holds for all the
conjugates of $\la_1$. Therefore, by \eqref{productformulaAr},
$f_k$ is divisible by all the factors $(1+e^{w\la_1})^{k-r}$.
Hence, it is also divisible by $(h+2)^{k-r}$. Since $h$ is
$W_0$-invariant, we obtain that $f_k$ is proportional to
$(h+2)^{k-r}\Delta$. This concludes the proof of the lemma.
\end{proof}
For the rest of the proof we need to divide $\overline{\ka_+}$ in a few parts. If $\{i_1,\dots,i_r\}=\{1,\dots,r\}$, we set
$$
\Lambda(i_1,\dots,i_r):=\{s\in \overline{\ka_+}\mid \bra\la_{i_1},s\ket\ \ge \dots \ge\ \bra\la_{i_r},s\ket\}.
$$
In fact only a few of these sets are enough to cover $\overline{\ka_+}$
\begin{lem}
The set $\Lambda(i_1,\dots,i_r)$ has a non empty intersection with $\ka_+$ if, and only if, for all $2\le j\le r$, there exists
$k<j$, such that $i_j=i_k\pm 1$.
\end{lem}
\begin{proof} Let $s=s_1\alpha_1+\dots +s_r\alpha_r$ be the
decomposition of $s\in \ka_+$ in the basis of the simple roots.
Assume that the hypothesis on $i_1,\dots,i_r$ of the lemma does
not hold. It means that for some $i$, we have either
$\bra\la_i,s\ket \ge \bra\la_{i+2},s\ket \ge \bra\la_{i+1},s\ket$,
or $\bra\la_{i+2},s\ket \ge \bra\la_i,s\ket \ge
\bra\la_{i+1},s\ket$. Assume that we are in the first case (the
second is similar). Then $s_i\ge s_{i+2}\ge s_{i+1}$. But on the
other hand, since $s\in \ka_+$, $\bra\alpha^\vee_{i+1},s\ket >0$.
This is absurd because
$\bra\alpha^\vee_{i+1},s\ket=s_{i+1}-\frac{1}{2}(s_i+s_{i+2})$.
\end{proof}
Now for the same
reason than in rank $2$ the equation $\nabla h(s)=\delta h(s)$ has a unique solution in $\overline{\ka_+}$. In the sequel
we will assume that it lies in $\Lambda(1,\dots,r)$. The proof works the same in the other cases.
\begin{lem}
\begin{enumerate}
\item For all $1\le i\le r$,
$$
e^{\bra\la_{i+1}-\la_i,s\ket}\asymp (1-\bra\alpha_1+\dots+\alpha_i,\delta\ket).
$$
\item For all $2\le i\le r$, when $\bra\la_i-\la_{i-1},s\ket\to 0$,
$$
\bra\la_i-\la_{i-1},s\ket\asymp \delta_i-\delta_{i-1}.
$$
\end{enumerate}
\end{lem}
\begin{proof} Let us prove the first claim by induction on $i\le
r$. The equation $\bra\alpha_1,\nabla h(s)\ket=\delta_1 h(s)$ may
be rewritten as
$$
\sum_\mu \sinh \bra\la_1+\mu,s\ket=\delta_1 (\sum_\mu \cosh \bra\la_1+\mu,s\ket+e^{\bra\la_2,s\ket}+\sum_\nu e^{\bra\nu,s\ket}),
$$
where the last sum is over weights $\nu$ such that $\bra\nu,s\ket \le \bra\la_2,s\ket$. Multiplying the left and right members
of the last equality by $e^{-\bra\la_1,s\ket}$ gives immediately $e^{\bra\la_2-\la_1,s\ket}\asymp 1-\delta_1$.
Let now $i\le r$. We write
$$
\bra\alpha_1+\dots +\alpha_i,\nabla h(s)\ket=\bra\alpha_1+\dots+\alpha_i,\delta\ket h(s).
$$
The exponential $e^{\la_1+\la_{i+1}-\la_i}$ appears in the right member of the last equality, whereas it does not in the left
member. Then we conclude by the same argument as before. The last statement is proved in a similar way, by using the equations
$\bra\alpha_{i+1}-\alpha_i,\nabla h(s)\ket=(\delta_{i+1}-\delta_i)h(s)$. This concludes the proof of the lemma.
\end{proof}
The end of the proof of Theorem \ref{theoestimateheatA} is now
completely similar to the rank $2$ case, if one avoids Section
\ref{lowerboundA2boundary}. Let us just notice that the quadratic
form appearing in Lemma \ref{quadraticformphase} is equal in
general to
$$
q_\delta(\teta)=\sum_{i=1}^r e^{\bra\la_{i+1}-\la_i,s\ket}\bra\la_{i+1}-\la_i,\teta\ket^2.
$$
We leave the other details of the proof to the reader.
\end{proof}

\section{Green's function estimate}
\subsection{statement of the result}
Green's function is defined for $x,y\in \kX$ by
$$
G(x,y|z)=\sum_{n\ge |x|} p^n(x,y)z^n,
$$
for all $z\in \C$ such that $|z|\le \frac{1}{\tilde{\rho}}$. We
set
$$
G(x,z):=G(0,x|z).
$$
We will give a sharp estimate of this function when $z$ is real
positive. As usual we always use the implicit notation $x\in
V_la(0)$ relating $x$ and $\la$.
\begin{theo}
\label{estimationgreen}
\begin{enumerate}
\item Let $z\in (0,\tilde{\rho}^{-1})$. Then
$$
G(x,z)\asymp \frac{1}{|\la|^{|R^+|+\frac{r-1}{2}}} e^{-\bra\la,s_0\ket}F_0(\la),
$$
where $s_0\in \overline{\ka_+}$ is uniquely determined by the
conditions: $h(s_0)=(\rho z)^{-1}$ and $\nabla h(s_0)$ is
proportional to $\la$.
\item We have
$$
G(x,\tilde{\rho}^{-1})\asymp \frac{1}{|\la|^{2|R^+|+r-2}}F_0(\la).
$$
\end{enumerate}
\end{theo}
\subsection{Proof}
\subsubsection{The case $z<\tilde{\rho}^{-1}$}
First we need some preliminary results. We set $g=\nabla \log h$, and for any $\delta \in \overline{\ka_+}$, we define
$s=s(\delta)$ as the unique point in $\overline{\ka_+}$ such that $g(s)=\delta$. We have the following
\begin{lem}
\label{diffg} The function $g$ is locally invertible. Moreover,
its differential $dg_s$ at any point $s\in \overline{\ka_+}$
satisfies $\bra u,dg_s(u)\ket >0$, for all $u\in \ka$.
\end{lem}
\begin{proof} We compute the differential of $g$ at some point
$s$
$$
dg_s(u)=\frac{1}{h^2(s)}\sum_{\la,\la'}[\frac{\bra\la,u\ket}{2}\la+\frac{\bra\la',u\ket}{2}\la'
                                             -\bra\la,u\ket\la'-\bra\la',u\ket\la] e^{\bra\la+\la',s\ket}.
$$
Thus
$$
\bra u,dg_s(u)\ket=\frac{1}{2h^2(s)}\sum_{\la,\la'}[\bra\la,u\ket-\bra\la',u\ket]^2e^{\bra\la+\la',s\ket}>0,
$$
for all $u\in \ka$. In particular $dg_s$ is invertible at each
point $s$. We conclude by using the local inversion theorem.
\end{proof}
For $t> |\la|$ we set $\delta_t=\frac{\la}{t}$, and
$s_t=g^{-1}(\delta_t)$. Now we define the function $\Psi$ in
$(|\la|,\infty)$ by
$$
\Psi(t)=t[\log h (s_t)-\bra\delta_t,s_t\ket+\log (\rho z)].
$$
We have
$$
\Psi'(t)=\log h(s_t)+\log (\rho z),
$$
and
$$
\Psi''(t)=-\bra\delta_t,dg_{\delta_t}^{-1}(\frac{\delta_t}{t})\ket.
$$
In particular, by Lemma \ref{diffg}, $\Psi''(t)<0$ for all $t>|\la|$. Thus $\Psi$ is strictly concave, and attains its maximum
at a unique point $t_0$, which satisfies the equation $h(s_{t_0})=(\rho z)^{-1}$. For simplify we set $s_0:=s_{t_0}$ and
$\delta_0=\delta_{t_0}$. Observe that they depend only on $\frac{\overline{x}}{|x|}$. We have
\begin{lem}
There exist constants $c>0$ and $C>0$ such that for all $x\in
\kX$,
$$
c\le |s_0|\le C \qquad \text{and}\qquad c\le |\delta_0| \le 1-c.
$$
\end{lem}
\begin{proof} The proof is straightforward. First $h(s_0)=(\rho
z)^{-1}$. Moreover $h$ is continuous, $h(0)<(\rho z)^{-1}$, and
$h(s)\to \infty$, when $|s|\to \infty$. Eventually $s\to 0$ when
$|\delta|\to 0$, and $|s|\to \infty$ when $|\delta|\to 1$.
\end{proof}
In the sequel it will be convenient to introduce also the function $\Phi$ defined for any $|\delta|<1$ by
$$
\Phi(\delta)= \log h(s)-\bra \delta,s \ket+\log (\rho z).
$$
We have
\begin{equation}
\label{diffPhi}
\nabla \Phi (\delta)=-s,
\end{equation}
and for all $u\in \ka$,
\begin{equation}
\label{diff2Phi}
d^2\Phi_{\delta}(u)=-\bra u,dg_\delta^{-1}(u)\ket,
\end{equation}
where $d^2\Phi_\delta$ is the second order differential of $\Phi$
at the point $\delta$. Writing the Taylor development of $\Phi$ at
order $2$, we get
$$
\Phi(\delta)=-\bra\delta,s_0\ket-\bra\delta-\delta_0,dg^{-1}_{\delta_0}(\delta-\delta_0)\ket+\kO(|\delta-\delta_0|^2).
$$
Thus there exists $\epsilon>0$, $c>0$ and $C>0$ such that
\begin{equation}
\label{estimationPhi}
-C|\delta-\delta_0|^2\le \Phi(\delta)+\bra\delta,s_0\ket \le -c|\delta-\delta_0|^2,
\end{equation}
for all $\delta$ such that $|\delta-\delta_0|\le\epsilon$. Taking smaller $\epsilon$ if necessary, we can also assume that there
exists a constant $c>0$ such that for $\delta$ in the preceding range, $c\le|\delta|\le 1-c$. We can now prove
\begin{lem}
\label{estimationprincipal}
We have the estimate
$$
\sum_{|\delta_n-\delta_0|\le \epsilon}p^n(0,x)z^n\asymp \frac{1}{|\la|^{|R^+|+\frac{r-1}{2}}} e^{-\bra\la,s_0\ket}F_0(\la).
$$
\end{lem}
\begin{proof} First by Theorem \ref{theoestimateheatA}, we see
that for all $n$ such that $|\delta_n-\delta_0|\le \epsilon$,
$$
p^n(0,x)z^n\asymp \frac{1}{|\la|^{|R^+|+\frac{r}{2}}} e^{\Psi(n)}F_0(\la).
$$
But $\Psi(n)=n\Phi(\delta_n)$. Thus by \eqref{estimationPhi}
$$
e^{-\bra\la,s_0\ket-Cn|\delta-\delta_0|^2}\le e^{\Psi(n)}\le e^{-\bra\la,s_0\ket-cn|\delta-\delta_0|^2},
$$
for all $n$ such that $|\delta_n-\delta_0|\le \epsilon$. Now we write
$$
\sum_{|\delta_n-\delta_0|\le \epsilon}e^{\Psi(n)}=\sum_{p=0}^{+\infty} \sum_{\frac{\epsilon}{2^{p+1}}
                                                 \le |\delta_n-\delta_0|
                         \le \frac{\epsilon}{2^p}}e^{\Psi(n)}.
$$
Next for all $p\ge 0$,
$$
\left|\left\{n\mid \frac{\epsilon}{2^{p+1}}\le|\delta_n-\delta_0|
                                           \le \frac{\epsilon}{2^p}\right\}\right| \asymp \frac{\epsilon|\la|}{2^p}.
$$
Thus
$$
\sum_{|\delta_n-\delta_0| \le \epsilon}e^{\Psi(n)}\le \text{const}\cdot |\la|^{\frac{1}{2}}e^{-\bra\la,s_0\ket}
 \sum_{p=0}^{+\infty} \frac{|\la|^{\frac{1}{2}}}{2^p}e^{-c|\la|/2^{2p}}.
$$
Moreover it is elementary to see that the last sum is bounded. Therefore we get as expected
$$
\sum_{|\delta_n-\delta_0|\le \epsilon}e^{\Psi(n)}\le C |\la|^{\frac{1}{2}}e^{-\bra\la,s_0\ket}.
$$
for some constant $C>0$. By the same argument we can prove the lower estimate
$$
\sum_{|\delta_n-\delta_0|\le \epsilon}e^{\Psi(n)}\ge c |\la|^{\frac{1}{2}}e^{-\bra\la,s_0\ket},
$$
for some constant $c>0$. This concludes the proof of the lemma.
\end{proof}
The proof is now almost finished. Since $\Psi$ is concave, there exists $c>0$ such that
$$
\Psi(n)\le -\bra\la,s_0\ket-c n,
$$
for all $n$ such that $|\delta_n-\delta_0|\ge \epsilon$. Thus
\begin{equation}
\label{etimationqueue1}
\sum_{|\delta_n-\delta_0|\ge \epsilon,\ |\delta_n|\ge \epsilon}e^{\Psi(n)}\le C|\la|e^{-\bra\la,s_0\ket-c|\la|}.
\end{equation}
In fact in the preceding sum we must assume that $n>|x|$, because the estimate of $p^{|x|}(0,x)$ is not contained in Theorem
\ref{theoestimateheatA}. But by the Harnack principle (given by the heat equation), there exists a constant $C>0$ such that for
all $x$,
$$
p^{|x|}(0,x)\le C p^{|x|+2}(0,x),
$$
which gives also a control of the term $p^{|x|}(0,x)$. Taking now smaller $\epsilon$ if necessary, we can assume that there
exists $c>0$ such that,
$$
\rho z h(s)\le 1-c,
$$
when $|\delta|\le \epsilon$. Then, by Theorem \ref{theoestimateheatA},
$$
\sum_{|\delta_n|\le \epsilon} e^{\Psi(n)}\le C \sum_{n \ge \frac{|\la|}{\epsilon}}(1-c)^n\le C(1-c)^{\frac{|\la|}{\epsilon}}.
$$
Taking again smaller $\epsilon$ if necessary, we get
\begin{equation}
\label{etimationqueue2}
\sum_{|\delta_n|\le \epsilon} e^{\Psi(n)}\le C e^{-\bra\la,s_0\ket-c|\la|}.
\end{equation}
Together with Lemma \ref{estimationprincipal}, and \eqref{etimationqueue1}, this proves Theorem \ref{estimationgreen} in
the case $z<\tilde{\rho}^{-1}$.

\subsubsection{The case $z=\tilde{\rho}^{-1}$}
First by Lemma \ref{diffg} $dg_0$ is invertible. Thus $|s|\asymp |\delta|$ near $0$. It implies with (\ref{diffPhi}) that there
exists a constant $c>0$ such that
$$
\Phi(\delta)\le -c|\delta|^2,
$$
for all $|\delta|<1$. Therefore for all $\epsilon>0$, we get a constant $c'>0$ such that
$$
\sum_{|\delta_n|\ge \epsilon} p^n(0,x)\tilde{\rho}^{-n}\le e^{-c'|x|}F_0(\la).
$$
Moreover, by \eqref{diffPhi} and \eqref{diff2Phi}, $\nabla \Phi (0)=0$ and $d^2\Phi_0$ is definite negative. Hence we get
$$
\Phi(\delta)\asymp -|\delta|^2,
$$
near $0$. It follows that
$$
\sum_{|\delta_n|\le \epsilon}p^n(0,x)\tilde{\rho}^{-n}\asymp F_0(\la)\int_{\frac{|x|}{\epsilon}}^{+\infty}
                                                             \frac{1}{t^m}e^{-\frac{|x|^2}{t}}dt,
$$
where $m=|R^+|+\frac{r}{2}$. Next we do the change of variable $t\to \frac{|x|^2}{t}$ and we find the desired estimate.
This concludes the proof of Theorem \ref{estimationgreen}.

\vspace{0.5cm}
\noindent J.-Ph. Anker and Br. Schapira \\
\noindent \textit{Universit\'e d'Orl\'eans, F\'ed\'eration Denis Poisson, Laboratoire MAPMO \\
B.P. 6759, 45067 Orl\'eans cedex 2, France.}\\
\noindent mails: jean-philippe.anker@univ-orleans.fr and bruno.schapira@univ-orleans.fr

\vspace{0.2cm}

\noindent Br. Schapira \\
\noindent \textit{Universit\'e Pierre et Marie Curie, Laboratoire de Probabilit\'es et Mod\`eles Al\'eatoires, \\
4 place Jussieu, F-75252 Paris cedex 05, France.}

\vspace{0.2cm}

\noindent B. Trojan \\
\noindent \textit{School of Mathematics and Statistics \\
Carslaw Building (F07), University of Sydney NSW 2006, Australia.} \\
\noindent mail: B.Trojan@maths.usyd.edu.au

\end{document}